\documentclass[12pt]{article}

\usepackage[T1]{fontenc}
\usepackage{amsthm, amssymb}
\usepackage{mathtools} %to declare the Paired Delimiters below
\usepackage{dsfont} %to define the characteristic function \ind below.
\usepackage{graphicx}
\usepackage{enumerate}
\usepackage{xcolor}

\newcommand{\R}{\mathbb{R}}

\newcommand{\Int}{\mathbb{Z}}
\newcommand{\Mod}[1]{\,(\textrm{mod}\; #1)}
\newcommand{\U}[1][t]{e^{i #1 \hbar\Delta/2}}

\DeclarePairedDelimiter\abs{\lvert}{\rvert}
\DeclarePairedDelimiter\norm{\lVert}{\rVert}

\DeclareMathOperator{\supp}{supp}

\DeclareMathOperator{\BigO}{\mathcal{O}}

\theoremstyle{plain}
  \newtheorem{theorem}{Theorem}
  
  \newtheorem{lemma}[theorem]{Lemma}
  
\theoremstyle{definition}  
  \newtheorem{definition}[theorem]{Definition}

\title{Convergence over fractals for the Schr\"odinger equation}

\author{Luc\`a, R. and Ponce-Vanegas, F.}

\date{}

\begin{document}
\maketitle

\begin{abstract}
We consider a fractal refinement of the Carleson problem for the Schr\"odinger equation, that is to identify the 
minimal regularity needed by the solutions to converge pointwise to their initial data almost everywhere with respect to the $\alpha$-Hausdorff 
measure ($\alpha$-a.e.). 
We extend to the fractal setting ($\alpha < n$) a recent counterexample of Bourgain \cite{Bourgain2016}, which is sharp in the Lebesque measure 
setting ($\alpha = n$).
In doing so we recover the necessary condition from \cite{zbMATH07036806} for pointwise 
convergence~$\alpha$-a.e. and we extend it to the range $n/2<\alpha  \leq (3n+1)/4$.
\end{abstract}

\section{Introduction}

A classic question related to solutions to the linear Schr\"odinger equation 
(here $\hbar = 1/(2\pi)$)
\begin{equation*}
\begin{dcases}
\partial_tu = i\frac{\hbar}{2}\Delta u \\
u(x,0) = f(x) \in H^s(\R^n),
\end{dcases}
\end{equation*}
is: what is the minimal regularity the initial datum must have so that the solution $u$ converge almost everywhere (a.e.) to $f$? 
%{\color{red} I reworded this part: whether 
%is which Sobolev regularity is needed in order to have convergence almost everywhere (a.e.) with respect to the Lebesgue measure to their initial data}. 
More precisely, which is the smallest $s\geq 0$ such that
\begin{equation}\label{LinConv}
\lim_{t \to 0}u(x,t) = f(x), \qquad \mbox{for a.e. $x \in \R^n$ and for all $f \in H^{s}(\R^n)$.} 
\end{equation}
This problem was introduced by Carleson in \cite{Carleson1980}, where he proved the validity of \eqref{LinConv} for $s \geq 1/4$ in dimension $n=1$.
Soon later Dahlberg and Kenig \cite{DahlbergKenig1982} proved this to be sharp. 
The considerably harder higher dimensional problem was subsequently studied by many authors 
\cite{Cowling1983, Carbery1985, Sjolin1987, Vega1988, Bourgain1992, MoyuaVargasVega1999, TaoVargas2000_1, TaoVargas2000_2, Tao2003, Lee2006, Bourgain2013, LucaRogers2019, DemeterGuo2016, LucaRogers2017, DuGuthLiZhang2018}. 

Recently, the problem has been settled, up to the endpoint, thanks to the contributions of Bourgain~\cite{Bourgain2016} (see \cite{Pierce2019} for a nice detailed exposition), who proved the necessity of  
$s \geq \frac{n}{2(n+1)}$, and of Du--Guth--Li~\cite{DuGuthLi2017} and of Du--Zhang \cite{DuZhang2019}, who proved the sufficiency of $s > \frac{n}{2(n+1)}$ in dimensions~$n = 2$ and~$n \geq 3$, respectively.
We mention that, besides Bourgain's counterexample, the necessity of $s \geq \frac{n}{2(n+1)}$ can be proved also by different counterexamples \cite{zbMATH07036806}.

In this paper we consider a fractal refinement of the Carleson problem. Given $\alpha \in (0, n]$, the goal is to identify the 
smallest $0 \leq s \leq n/2$ such that
\begin{equation}\label{FractalLinConv}
\lim_{t \to 0}u(x,t) = f(x), \qquad \mbox{for $\alpha$-a.e. $x \in \R^n$ and for all $f \in H^{s}(\R^n)$,} 
\end{equation}
where $\alpha$-a.e. means almost everywhere with respect to the $\alpha$-dimensional Hausdorff measure.

%We explain how to do this in the next section. 
This fractal refinement of the Carleson problem was introduced in \cite{SjogrenSjolin1989}. In 
\cite{BarceloBennettCarberyRogers2011}, the authors gave a complete solutions for $\alpha \in [0,n/2]$, proving that 
$s > (n-\alpha)/2$  is necessary and sufficient for \eqref{FractalLinConv} to hold. 
The necessity of this condition depends on the Sobolev space framework, since for smaller $s$
there exist initial data in $H^s(\mathbb R^n)$ that are not well defined on sets of 
dimension $\alpha$; see~\cite{Zubrinic2002}. 
On the other hand, for $s > (n-\alpha)/2$ one can make sense of the initial data and of the relative solution
$\alpha$-a.e.; we refer to the proof of Theorem \ref{thm:Example} for details.  
When $\alpha \in (n/2, n]$, Du and Zhang \cite{DuZhang2019} proved the best known sufficient condition for \eqref{FractalLinConv} to hold:
\begin{equation}\label{FractCond1}
s > \frac{n}{2(n+1)} (n + 1 - \alpha).
\end{equation}
As mentioned, this is optimal (up to the endpoint) when $\alpha =n$, 
but it is not clear yet whether this is optimal for $\alpha$ strictly smaller. It is worth mentioning that
\eqref{FractCond1} is necessary for the $\alpha$-a.e. pointwise convergence in the periodic setting \cite{ecei2020}, 
however in this setting it is still unknown if it is sufficient (not even for $\alpha =n$).
 
In \cite{zbMATH07036806} it was proved that for 
$(3n+1)/4 \leq  \alpha \leq n$ the condition
\begin{equation}\label{INTVAL}
s > \frac{n}{2(n+1)}+\frac{n-1}{2(n+1)}(n-\alpha) \, ,
\end{equation}
is necessary for \eqref{FractalLinConv} to hold. Here we extend this result to the full range~$n/2 < \alpha \leq n$ (recall that for smaller $\alpha$ the problem has been solved in \cite{BarceloBennettCarberyRogers2011}); thus the result is new for $n/2 < \alpha \leq (3n+1)/4$. To prove this result, we use
a modification of the Bourgain counterexample rather than the counterexample in \cite{zbMATH07036806}. 
We consider this fact of independent interest. The possibility of adapting the Bourgain counterexample to the 
fractal measure setting was also suggested by Lillian Pierce in \cite{Pierce2019}.

\begin{theorem}\label{thm:Main}
Let $n\ge 2$ and $n/2 < \alpha \leq n $. Then
for every 
\begin{equation}\label{Strict}
s'<s := \frac{n}{2(n+1)}+\frac{n-1}{2(n+1)}(n-\alpha)  
\end{equation}
there exists a function $f\in H^{s'}(\R^n)$ such that
\begin{equation}\label{Eq:Stat}
\limsup_{t\to 0^+}\abs{\U f(x)} = \infty
\end{equation}
for $x$ in a set of Hausdorff dimension $\geq \alpha$.  
\end{theorem}

For $\alpha \in (n/2, n)$ we can in fact immediately improve the statement, saying that \eqref{Eq:Stat} occurs on a set with $\alpha$-Hausdorff measure $= \infty$. 
This is because in \eqref{Strict} we have a strict inequality. Thus, given $\alpha' > \alpha$ and sufficiently close to $\alpha$ in such a way that
$$
s'<s := \frac{n}{2(n+1)}+\frac{n-1}{2(n+1)}(n-\alpha'),  
$$
we would in fact prove that \eqref{Eq:Stat} occurs on a set of dimension $\geq \alpha'$.
When~$\alpha = n$ we can not self-improve the statement, however we know by \cite{zbMATH07036806} that \eqref{Eq:Stat} holds 
on a set of strictly positive Lebesgue measure.

A consequence of Theorem \ref{thm:Main} is the necessity of the condition   
$$
s \geq  \frac{n}{2(n+1)}+\frac{n-1}{2(n+1)}(n-\alpha) \, 
$$
for the validity of the maximal estimate 
\begin{equation}\label{WIN1}
\int_{B_R} \sup_{t \in (0,1)} |e^{it\Delta} f(x)|^2 d \mu(x)  \lesssim C_\mu R^{2s} \|f\|^2_{2} \,,
\end{equation}
where $B_R\subset\R^n$ is a ball of radius $R>1$, and $\mu$ is an $\alpha$-dimensional measure on $B_R \subset \R^n$, \textit{i.e.} a positive Borel measure that satisfies
$$\mu(B_r(x)) \lesssim C_\mu  r^{\alpha}, $$ for all balls with center $x$ and radius $r >0$.
One may see \eqref{WIN1} as the weighted $L^2$ inequality  
\begin{equation}\label{WIN2}
 \int_{B_R} \sup_{t \in (0,1)} |\widehat{g d \sigma}(x)|^2  d \mu(x)  \lesssim C_\mu R^{2s} \|g\|^2_{L^{2}(S)} ,
\end{equation}
where $S$ is a bounded hypersurface in $\R^d := \R^{n+1}$ with non zero gaussian curvature (for instance, a portion of 
the paraboloid in the case of 
\eqref{WIN1}) and $d\sigma$ is the measure induced on $S$ by the Lebesgue measure. A closely related family of weighted $L^2$ estimates is 
\begin{equation}\label{WIN3}
 \int_{B_1}  |\widehat{g d \sigma}(Rx)|^2  d \mu(x) \lesssim C_\mu R^{-\gamma} \|g\|^2_{L^{2}(S)},  
\end{equation}
where $B_1$ is now a ball in $\R^d$ of radius $1$, 
and $\mu$ is an $\alpha$-dimensional measure on $B_1 \subset \R^d$.  The problem here is to identify the largest $\gamma$ such 
that~\eqref{WIN3} holds. 
Interestingly, these problems are very sensitive to the arithmetical structure of the hypersurface $S$.
 For instance, the known necessary conditions are different
for the sphere and the paraboloid; see \cite{Iosevich2007, BarceloBennett2007, LucaRogers2019, Du2020, ponce2018, zbMATH07160159}.

%{\color{red} I am not sure if it is possible with Bourgain's counterexample. I know that with your example [19] it is indeed possible, but a technical issue with the totient function makes it hard to sustain the claim with Bourgain's example, at least when $n\ge 3$; I checked in Pierce's paper, sec. 5.1(II) and sec. 7, and she seems not to get rid of this problem; the volume of divergence is smaller $\sim R^{-\varepsilon}_k$ as $k$ goes to infinity. I don't know how Bourgain justified it rigorously. Quite probably I ignore something, so I'd be happy if you clarify this issue to me. Thank you.} We will omit the details of this last fact.

%\begin{theorem}\label{thm:Main}
%If $n\ge 2$ and $\frac{n}{2(n+1)}\le s<\frac{n}{4}$, then there exists a function $f\in H^{s'}(\R^n)$, for every $s'<s$, such that
%\begin{equation}
%\limsup_{t\to 0^+}\abs{\U f(x)} = \infty
%\end{equation}
%in a set of Hausdorff dimension 
%\begin{equation}
%\alpha = n+\frac{n}{n-1}-2\frac{n+1}{n-1}s.
%\end{equation}
%\end{theorem}

\subsection*{Notations}

\begin{itemize}
\item $e(z) = e^{iz}$.
\item If $A\subset\R^n$, then $\abs{A}$ is its Lebesgue measure, and if $A$ is a discrete set, then $\abs{A}$ is the cardinality. For example, if $I = [a,b]\subset \mathbb{Z}$ denotes the interval of integers $a\le k\le b$, then $|I|$ is the length of the interval.
\item If $I = [a,b]\subset \mathbb{Z}$, for $a,b\in\mathbb{R}$, denotes an interval of integers, then we write $L(I) := \min_{k\in I} k$ and $R(I) := \max_{k\in I} k$.
\item $B_r(x)\subset\R^n$ is a ball of radius $r$ and center $x$---the center is usually omitted. $Q(x,l)\subset\R^n$ is a cube with side-length $l$ and center $x$. 
\item If $x\lesssim y$, then $x\le Cy$ for some constant $C>0$, and similarly for $x\gtrsim y$; if $x\simeq y$ then $x\lesssim y\lesssim x$. If $x\ll y$ then $x\le cy$, where $c$ is a sufficiently small constant, and similarly for $x\gg y$. 
\item $\limsup_{k\to \infty} F_k := \bigcap_{N\ge 1}\bigcup_{k\ge N} F_k$.
\item Hausdorff dimension of a set: for $0<\alpha\le n$ and $\delta>0$ we define the outer measure
\begin{equation*}
\mathcal{H}^\alpha_\delta(F) := \inf\{\sum_{B_r\in\mathcal{B}}r^\alpha\mid F\subset\bigcup_{B_r\in\mathcal{B}}B_r\textrm{ and } r<\delta\};
\end{equation*}
we do not exclude the case $\delta = \infty$. The $\alpha$-dimensional Hausdorff measure of a set $F$ is $\mathcal{H}^\alpha(F):=  \lim_{\delta\to 0}\mathcal{H}^\alpha_\delta(F)$. The Hausdorff dimension of a set $F$ is $\sup\{\alpha\mid \mathcal{H}^\alpha(F) > 0\}$.
\end{itemize}

\subsection*{Acknowledgments}

This research is funded by the Basque Government through the BERC 2018-2021 program, and by the Spanish State Research Agency through BCAM Severo Ochoa excellence accreditation SEV-2017-0718 and by the IHAIP project PGC2018-094528-B-I00. Additionally, the first author is supported by Ikerbasque and 
the second author is supported by the ERCEA Advanced Grant 2014 669689-HADE.

\section{Preliminaries}

We recall some classic estimates about exponential sums that we will use repeatedly in the rest of the paper.

%Hereafter $I$ will be an 
%an integer interval, which means that 
%its boundary points are natural numbers, %We also allow the right extreme to be~$\infty$. 
%namely $I = [K_1, K_2]$ with $K_1, K_2 \in \N$.
%We denote 
%by $|I|$ the length of~$I$ and we denote by $L(I) = K_1$ its left boundary and by $R(I)=K_2$ its right boundary. In general, we will use the short notation $k \in I$ for $k \in \mathbb{N} \cap I $. 

%We first need the following 

We recall first a classical result about Gauss quadratic sums, whose proof can be consulted in Lemma~3.1 of \cite{Pierce2019}.

\begin{lemma}[Gauss quadratic sums] \label{thm:GQS}
If $a,b,q \in \mathbb{Z}$ satisfy the conditions $(a, q) = 1$ and
\begin{equation}\label{eq:thm:fractions_GS}
\begin{cases}
b \in \mathbb{Z}  & \textrm{when } q \textrm{ is an odd number,} \\
b \textrm{ is even} & \textrm{when } q\equiv 0 \Mod 4, \\
b \textrm{ is odd} & \textrm{when } q\equiv 2 \Mod 4,
\end{cases}
\end{equation}
then for the quadratic phase
\begin{equation}\label{eq:thm:quadratic_phase}
f(r)  := \frac{a}{q} r^2 +  \frac{b}{q} r
\end{equation} 
it holds that
\begin{equation}\label{QGSFinal}
\left| \sum_{r =0}^{q-1} e(2\pi f(r) ) \right|  = c_q \sqrt{q}, % \qquad  |I|=q-1,
\end{equation}
where $c_q=1$ when $q$ is odd, and $c_q =\sqrt{2}$ when $q$ is even. 
\end{lemma}

The following estimate due to Weyl will be useful to handle incomplete Gauss sums. 

\begin{lemma}%[Van Der Corput]
\label{WeylBound}
Let $I$ be an integer interval. If $a,b,q \in \mathbb{Z}$  satisfy the conditions $(a, q) = 1$ and \eqref{eq:thm:fractions_GS}, then for the quadratic phase $f$ in \eqref{eq:thm:quadratic_phase} it holds that 
 \begin{equation}\label{eq:FullWeylBound}
\left| \sum_{k \in I } e(2 \pi  f(k)) \right| =   C\frac{|I|}{\sqrt{q}}  + \mathcal{O}( \sqrt{q\ln q} ),
\end{equation}
where $\frac{1}{2}\le C\le \sqrt{2}$.
%In particular, if $|I| \leq q$, 
%\begin{equation}\label{eq:WeylBound}
%\left| \sum_{k \in I } e^{2 \pi i f(k)} \right| \lesssim  \sqrt{q \ln q} 
%\end{equation}
\end{lemma}

%{\color{red} I made some modifications in the proof to avoid technical problems, especially the case $M_1 = M_2$. The code of your proof is inside comments. In the Lemma I replaced $\simeq$ by $=$.}

\begin{proof}
We can assume that $L(I) := \min_{k\in I} k = 0$. In fact, 
\begin{equation*}
\sum_{k = L(I)}^{R(I)} e(2\pi(\frac{a}{q}k^2+\frac{b}{q}k)) = e(2\pi(\frac{a}{q}L(I)^2 + \frac{b}{q}L(I)))\sum_{k = 0}^{|I|-1} e(2\pi(\frac{a}{q}k^2+\frac{b+2aL(I)}{q}k)),
\end{equation*} 
and the absolute value at both sides is the same; we observe that the parity of $b$ and $b+2aL(I)$ is preserved.

If $\abs{I}<q$, then
\begin{equation}\label{RemainderJ}
\left| \sum_{k \in I } e(2 \pi  f(k)) \right| \le C\sqrt{q\ln q};
\end{equation}
for the proof we refer to Lemma~3.2 of \cite{Pierce2019}.

If $\abs{I}\ge q$, then we can sum in blocks of length $q$. Let $M$ be the largest integer that satisfies $Mq\le \abs{I}$, \textit{i.e.} $Mq\le \abs{I}< (M+1)q$, then
\begin{align*}
I &= [0, Mq-1]\cup J \\
&= \Big(\bigcup_{m=0}^{M-1}[mq, mq + q-1]\Big)\cup J,
\end{align*}
where $\abs{J}<q$. The sum over each block $[mq,mq+q-1]$ is a Guass quadratic sum, and we arrive to
\begin{equation*}
\sum_{k \in I } e(2 \pi f(k)):=  M \sum_{r=0}^{q-1}  e(2 \pi f(r)) + \sum_{k\in J}  e(2 \pi f(k)).
\end{equation*}
By our election of $M$ we have $M = C\abs{I}/q$, for $\frac{1}{2}<C\le 1$, and by \eqref{RemainderJ} we have
\begin{equation*}
\left|\sum_{k \in I } e(2 \pi f(k))\right| = C\frac{\abs{I}}{q}\left| \sum_{r=0}^{q-1}  e(2 \pi f(r))\right| + \BigO(\sqrt{q\ln q}).
\end{equation*}
Finally, we apply Lemma~\ref{thm:GQS} to get \eqref{eq:FullWeylBound}.
\end{proof}

To deal with perturbations of quadratic sums, we will use the following Lemma, which is consequence of Abel's summation formula; see Lemma 2.3 of \cite{ecei2020}.
\begin{lemma}%[Abel]
\label{DirichletSumLemma}
Let $I$ be an integer interval. Let $ a_k \geq 0$ be a sequence of real numbers and $b_k$ be sequences of complex numbers such that 
\begin{enumerate}
\item $a_{k+1} \leq a_k,$ 
\item
$\left| \sum_{k \in I' } b_k  \right| \leq \mathcal{C} , \quad  \mbox{for every interval $I' \subseteq I$}$.
\end{enumerate}
Then, 
\begin{equation}\label{DirichletTest}
\left| \sum_{k \in I'} a_k b_k \right| \leq  \mathcal{C} a_{L(I')} , \quad  \mbox{for every interval $I' \subseteq I$}.
\end{equation}
If (1) is replaced with $a_{k+1} \geq  a_k$, then 
$$
\left| \sum_{k \in I'} a_k b_k \right| \leq  \mathcal{C} a_{R(I')}, \quad  \mbox{for every interval $I' \subseteq I$}.
$$
\end{lemma}
%\begin{proof}
%Let $I' = [K_1, K_2]$.
%We only consider the case $a_{k+1} \leq a_k$, since the case $a_{k+1} \geq a_k$ can be handled by changing variables
%$k \leftrightarrow K_1 + K_2 - k$. 
%The inequality~\eqref{DirichletTest} easily follows by the summation by parts formula
%\begin{equation}
%\sum_{k=K_1}^{K_2} a_k b_k = a_{K_2} B_{K_2} + \sum_{k=K_1}^{K_2-1} B_k(a_k - a_{k+1}), 
%\quad \text{ where } \quad 
%B_{k} = \sum_{\ell=K_1}^{k}b_\ell.
%\end{equation}
% By \textit{(2)}, we have    $|B_{k}| \leq \mathcal{C}$ for every $k$, so  
%\begin{equation}
%\begin{split}
%\left| \sum_{k=K_1}^{K_2} a_k b_k  \right| &  \leq \mathcal{C}\, |a_{K_2}| + \mathcal{C}\,  \sum_{k=K_1}^{K_2-1} |a_k - a_{k+1}|  = \mathcal{C} \left(  a_{K_2} + \sum_{k=K_1}^{K_2-1} (a_k - a_{k+1}) \right) \\
%& = \mathcal{C} \left(  a_{K_2} + a_{K_1} - a_{K_2}  \right) = \mathcal{C}\, a_{K_1}.
%\end{split}
%\end{equation}  
%\end{proof}

\section{The main lower bound}

The initial data we consider are modifications of the Bourgain's counterexample in \cite{Bourgain2016}. Let $\varphi$ be a smooth positive function such that $\supp\hat{\varphi}\subset B_1(0)$ and $\varphi(0) = 1$.  We define the function
\begin{equation}\label{eq:Initial_Datum}
f_D(x) := f_1(x_1) \tilde{f}(\tilde{x}) 
\end{equation}
where 
$$
f_1(x_1) = 
 e(2\pi Rx_1)\varphi(R^\frac{1}{2} x_1), \qquad 
\tilde{f}(\tilde{x}):=  \prod_{j=2}^n\varphi(x_j)\Big(\sum_{\frac{R}{2D} < l_j < \frac{R}{D}}e(2\pi Dl_jx_j)\Big)
$$
where
$l = (l_2,\ldots,l_n)\in \Int^{n-1}$ and $x = (x_1, \tilde{x})\in \R\times \R^{n-1}$. For now we set $D$ as a free parameter, and we will choose its value later as a suitable power of $R$.

We need the following definition before investigating the divergence set of $\U f_D$; compare with \eqref{eq:thm:fractions_GS}.

\begin{definition}[Admissible fractions]\label{def:Admissible_fractions}
Let $p_1, \ldots, p_n, q \in \mathbb{Z}$. A point $(p_1/q,\ldots, p_n/q)$ is an admissible fraction if $(p_1,q)=1$ and if
\begin{equation}\label{eq:def:admissible}
\begin{cases}
(p_2,\ldots,p_n)\in \Int^{n-1} & \textrm{when } q \textrm{ is an odd number,} \\
p_j \textrm{ are even} & \textrm{when } q\equiv 0 \Mod 4, \\
p_j \textrm{ are odd} & \textrm{when } q\equiv 2 \Mod 4.
\end{cases}
\end{equation}
\end{definition}

\begin{figure}
\centering
\includegraphics[scale=1]{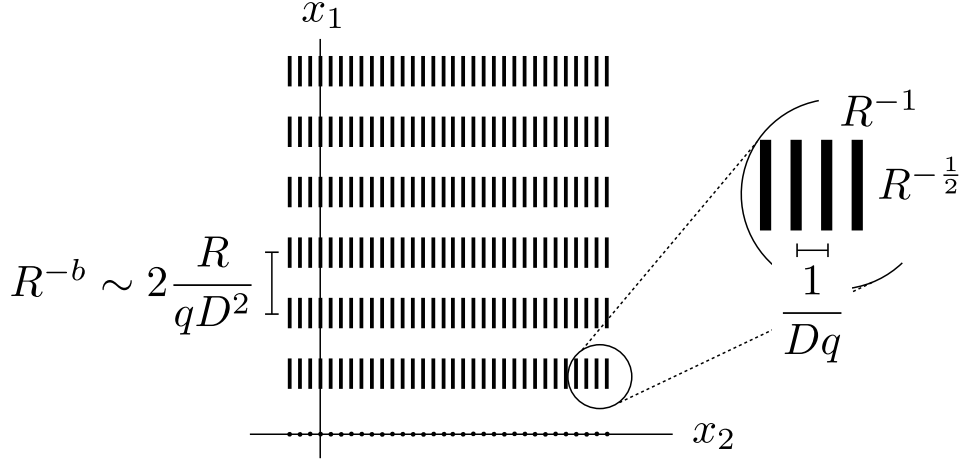}
\caption{Set $E_{q,D}$ in Theorem~\ref{thm:Supremum_Fq}. Some slabs may disappear to satisfy the conditions of admissibility.}\label{fig:Set_Eq}
\end{figure}

\begin{theorem}\label{thm:Supremum_Fq}
Let $c\ll 1$ and let $q>0$ be an integer such that $\frac{R}{Dq} \gg \sqrt{\ln{q}}$. If $f$ is the initial datum \eqref{eq:Initial_Datum}, then 
\begin{equation}\label{eq:thm:supremum_Fq}
\frac{\abs{\U f_D(x)}}{\norm{f_D}_{L^2}}\gtrsim R^\frac{1}{4}\Big(\frac{R}{Dq}\Big)^\frac{n-1}{2}
\end{equation}
for $(x,t)$ such that $0<t = 2 p_1/(D^2q)\ll 1/R$ and
\begin{equation}\label{eq:thm:Set_Box}
x \in E_{q,D} \cap [0,c]^n,  
\end{equation}
 where $E_{q,D}$ is the set of points
\begin{equation}\label{eq:thm:Set_F_q}
x_1 \in 2\frac{p_1R}{qD^2} + [-c R^{-\frac{1}{2}}, c R^{-\frac{1}{2}}] \quad\mbox{and}\quad  x_j \in  \frac{p_j}{Dq} + [-c R^{-1}, c R^{-1}],  \ \ 2\le j\le n;
\end{equation} 
here $(p_1/q,\ldots, p_n/q)$ is an \textit{admissible fraction} in the sense of Definition \ref{def:Admissible_fractions}; see Fig.~\ref{fig:Set_Eq}.
\end{theorem}

%{\color{red} I think it can happen that there is not $(x,t)$ that satisfies the conditions of the Theorem, but that is not a problem because our restrictions guarantee that the divergence set is not void.}
%
%Before proving the theorem, we make a comment on how the the choice~\eqref{eq:thm:Set_Box}
%imposes restrictions on $p_1$ and $p_j$. It implies indeed that~$|p_1| \leq  c \frac{q D^2}{2R}$ and that~$|p_j| \leq  c D q$. Note that
% $|p_1| \leq c \frac{q D^2}{2R}$ and the choice of $t = 2 p_1/(D^2q)$ are consistent with the assumption $|t| \leq \frac{c}{R}$.  

\begin{proof}
If $\hat{f}$ is an integrable functions, the solution of the Schr\"odinger equation with initial datum $f$ can be represented as
\begin{equation*}
e^{it\hbar\Delta/2}f(x) = \int \hat{f}(\xi)e(-\pi t\abs{\xi}^2 + 2\pi x\cdot \xi)\,d\xi.
\end{equation*}

We want to compute the modulus of $\U f_D(x)$ in the region $\abs{x}<c$ and $0< t<c/R$. We note that
$$
| \U f_D(x) | = | \U f_1(x_1) | \, | \U \tilde{f}(\tilde{x}) | \,.
$$
A direct computation shows that for 
$|t| \leq c/R$ and 
\begin{equation}\label{P3}
x_1 \in  tR + [-c R^{-\frac{1}{2}}, c R^{-\frac{1}{2}}]
\end{equation}
we have 
\begin{equation}\label{eq:f_DPreq}
| \U f_1(x_1) | \simeq \abs{\varphi(R^\frac{1}{2}(x_1-tR))} \simeq 1 
\end{equation}
Again, a direct computation gives ($\tilde{x} \in \R^{n-1}$)
\begin{equation}\label{eq:f_D}
\begin{split}
\U \tilde{f}(\tilde{x}) &=  \prod_{j=2}^n\int  \hat{\varphi}(\xi_j)e(-\pi t \xi_j^2 + 2\pi x_j \xi_j ) \\
&\hspace*{2cm} \sum_{\frac{R}{2D} < l_j< \frac{R}{D}} 
e(-\pi t\abs{Dl_j}^2 + 2\pi Dl_j(x_j-t\xi_j))\,d\xi_j.
\end{split}
\end{equation}
To estimate the absolute value of this product, we recall our hypotheses \eqref{eq:thm:Set_F_q}: $x_j = p_j/(Dq) + \varepsilon_j$, for $\abs{\varepsilon_j}<c/R$. We split each factor in \eqref{eq:f_D} into the main term
\begin{equation}\label{eq:main_GS}
F_{\textrm{main}}(t,p_j/q) := \U \varphi(x_j)\sum_{\frac{R}{2D} < l_j< \frac{R}{D}} 
e(-\pi t\abs{Dl_j}^2 + 2\pi l_j\frac{p_j}{q})
\end{equation}
and the perturbation
\begin{multline}\label{eq:perturbation_GS}
F_{\textrm{per}}(t,x_j) := \int  \hat{\varphi}(\xi_j)e(-\pi t \xi_j^2 + 2\pi x_j \xi_j ) \\
\sum_{\frac{R}{2D} < l_j< \frac{R}{D}} 
e(-\pi t\abs{Dl_j}^2 + 2\pi l_j\frac{p_j}{q})(1-e(2\pi Dl_j(\varepsilon_j-t\xi_j))\,d\xi_j.
\end{multline}

By hypothesis $t = 2p_1/(D^2q)$, so we can exploit Lemma~\ref{WeylBound} and the condition $R/(Dq)\gg \sqrt{\ln q}$ to estimate the main contribution \eqref{eq:main_GS} as
\begin{align}
\abs{F_{\textrm{main}}} &\simeq \abs{\sum_{\frac{R}{2D} < l_j< \frac{R}{D}} 
e\big(-2\pi(\frac{p_1}{q}l_j^2 - \frac{p_j}{q}l_j)\big)} \notag \\
&\simeq \frac{R}{D\sqrt{q}};\label{eq:main_GS_estimate}
\end{align}
we used $\abs{\U \varphi(x_j)}\simeq 1$.

We claim that the perturbation term \eqref{eq:perturbation_GS} satisfies $\abs{F_{\textrm{per}}}\ll R/(D\sqrt{q})$, which, together with \eqref{eq:f_D} and  \eqref{eq:main_GS_estimate}, leads to
\begin{equation}\label{I-Bound}
\abs{\U \tilde{f}(\tilde{x})} \simeq \Big(\frac{R}{D\sqrt{q}}\Big)^{n-1}.
\end{equation}
Then, we multiply by \eqref{eq:f_DPreq} to reach
\begin{equation}\label{Plug1}
\abs{\U f_D(x)} \simeq \Big(\frac{R}{D \sqrt{q}}\Big)^{n-1}.
\end{equation}
%We will show that 
%if we furthermore restrict to
% $$x_1 \in tR  + [-cR^{-\frac{1}{2}}, cR^{-\frac{1}{2}}] = \frac{2Rp_1}{D^2q}  + [-cR^{-\frac{1}{2}}, cR^{-\frac{1}{2}}],$$
%which is equivalent to the restrictions in \eqref{P3} and  \eqref{eq:thm:supremum_Fq}, and to $x_j$ as in \eqref{eq:thm:Set_F_q} we have
%\begin{equation}\label{I-Bound}
%|\eqref{eq:f_D}| 
%\simeq \Big(\frac{R}{D \sqrt{q}}\Big)^{n-1}
%\end{equation}
%By \eqref{eq:f_DPreq}--\eqref{I-Bound} we obtain
%Now \eqref{Plug1} gives \eqref{eq:thm:supremum_Fq}, 
 %after dividing by 
Finally, we divide \eqref{Plug1} by $\norm{f_D}_2 \simeq R^{-\frac{1}{4}}(R/D)^\frac{n-1}{2}$ to obtain \eqref{eq:thm:supremum_Fq}, and so the statement of the Theorem follows up to the claim $\abs{F_{\textrm{per}}}\ll R/(D\sqrt{q})$.

To prove the upper bound $\abs{F_{\textrm{per}}}\ll R/(D\sqrt{q})$, where $F_{\textrm{per}}$ was defined in \eqref{eq:perturbation_GS}, we begin with
\begin{align*}
\abs{F_{\textrm{per}}} &\lesssim \sup_{\abs{\xi_j}\le 1}\abs{\sum_{\frac{R}{2D} < l_j< \frac{R}{D}} 
e(-2\pi \frac{p_1}{q}l_j^2 + 2\pi l_j\frac{p_j}{q})(1-e(2\pi Dl_j(\varepsilon_j-t\xi_j))} \\
&= \sup_{\abs{\xi_j}\le 1}\Bigl|  \sum_{\frac{R}{2D} < l_j<\frac{R}{D}} 
e(- 2 \pi  \frac{p_1}{q} l_j^2 + 2\pi \frac{p_j}{q} l_j  ) \phi_{\varepsilon,t,\xi_j}(l_j) \Bigr|
\end{align*}
where $(p_1,q) = 1$, and $\phi_{\varepsilon_j,t,\xi_j}(l_j) = 1 - e(2\pi  D l_j (\varepsilon_j - t \xi_j)  )$.

By the triangle inequality, it suffices to prove
\begin{equation}\label{ISFinal}
\Bigl|  \sum_{\frac{R}{2D} < l_j<\frac{R}{D}} 
e(- 2 \pi  \frac{p_1}{q} l_j^2 + 2\pi \frac{p_j}{q} l_j ) \phi^i_{(\cdot)}(l_j) \Bigr| \lesssim c \frac{R}{D \sqrt{q}}  , \qquad c \ll 1 \qquad i=1,2.
\end{equation}
where
$$
\phi^1_{(\cdot)}(l_j) := 
 1 - \cos(2\pi  D l_j (\varepsilon_j - t \xi_j) ) 
\enspace\mbox{and}\enspace
\phi^2_{(\cdot)}(l_j) := 
%\frac{\xi_j}{| \xi_j|} \sin ( 2\pi  D l_j (\varepsilon_j - t \xi_j) ) \, .
\abs{\sin ( 2\pi  D l_j (\varepsilon_j - t \xi_j) )} \, .
$$
Again by Lemma \ref{WeylBound}, and using $R/(Dq) \gg \sqrt{\ln{q}}$, we have that
\begin{align}
\Bigl| \sum_{l_j \in I }   e(- 2 \pi  \frac{p_1}{q} l_j^2 + 2\pi \frac{p_j}{q} l_j ) 
 \Bigr| 
& \lesssim \frac{|I|}{\sqrt{q}}  +   \sqrt{q\ln q}
 \\ \nonumber
 & \lesssim \frac{R}{D \sqrt{q}} +   \sqrt{q\ln q} \lesssim \frac{R}{D \sqrt{q}} , 
 \quad \forall \, |I| \leq \frac{R}{D} \, .
\end{align}
On the other hand, the functions $\phi^i_{(\cdot)}(l_j)$ are real valued, positive, increasing in $l_j$, and satisfy
$$
\phi^i_{(\cdot)}(l_j) \lesssim
|D| |l_j| |\varepsilon_j - t \xi_j| \lesssim c, \qquad c \ll1;
$$
recall that 
$|l_j| \leq \frac{R}{D}$, $|\varepsilon_j| \leq \frac{c}{R}$, $|t| \leq \frac{c}{R}$ and $\abs{\xi_j} \leq 1$ in the support of $\hat \phi_j$. 
Thus~\eqref{ISFinal} follows by the second part of Lemma \ref{DirichletSumLemma} taking   
$$
a_{l_j} =  \phi^i_{(\cdot)}(l_j) 
\qquad \mbox{and}\qquad
b_{l_j} = e(- 2 \pi  \frac{p_1}{q} l_j^2 + 2\pi \frac{p_j}{q} l_j ),
$$
and the proof is concluded.

\end{proof}

\section{Construction of the Examples}

According to Theorem~\ref{thm:Supremum_Fq}, the function $\sup_{0<t<1}\abs{\U f_D}$ is large in the set
\begin{equation*}
\bigcup_{1\le q \le Q} (E_{q,D} \cap [0,c]^n) \subset\R^n, \quad 0 < c  \ll 1,
\end{equation*}
as long as $\frac{R}{DQ}\gg \sqrt{\ln Q}$. 

To cover the largest possible area, we should ensure that the collection of sets $E_{q,D}$, for $1\le q\le Q$, is essentially pairwise disjoint. In a unit cell $[0,1/D]^{n-1}$, the number of fractions $(p_2/(Dq),\ldots,p_n/(Dq))$, for $1\le q\le Q$, is $\simeq Q^n$, and if we think of the fractions as if they were uniformly distributed, then the average distance between them is $\simeq 1/(Q^\frac{n}{n-1}D)$, so we impose the restriction
\begin{equation}\label{eq:Condition_a}
R^{-a} := \frac{1}{Q^\frac{n}{n-1}D} \ge R^{-1}.
\end{equation}
We remark that $R/(DQ) = Q^\frac{1}{n-1}R^{1-a}$, so the condition $R/(Dq)\gg \sqrt{\ln q}$, for $1\le q\le Q$, is easily satisfied. 

The slabs that form $E_{q,D}$ have dimensions $cR^{-\frac{1}{2}}\times cR^{-1}\times\cdots\times cR^{-1}$, for $c\ll 1$, and they do not overlap in the $\tilde{x}$-space because $R/(DQ)\gg 1$, however they may overlap in the $x_1$ direction. To exploit the whole area of the slabs, we impose the new restriction
\begin{equation}\label{eq:Condition_b}
R^{-b} := \frac{R}{QD^2}\ge R^{-\frac{1}{2}};
\end{equation}
see Figure~\ref{fig:Set_Eq}.

The conditions \eqref{eq:Condition_a} and \eqref{eq:Condition_b} allow us to solve for  $Q$ and $D$ as
\begin{equation}\label{eq:D_Q}
D = R^{(n-(n-1)a+nb)/(n+1)}\qquad\mbox{and}\qquad Q = R^{\frac{n-1}{n+1}(2a-b-1)}.
\end{equation}
Since $Q\ge 1$, then we have to be sure that $2a\ge 1+b$, so we can write our conditions as
\begin{equation}\label{eq:conditions_a_b}
0<a\le 1, \quad 0<b\le\frac{1}{2} \quad\mbox{and}\quad 2a\ge 1+b;
\end{equation}
in particular, $a\ge\frac{1}{2}$.

\begin{center}
\includegraphics[scale=0.6]{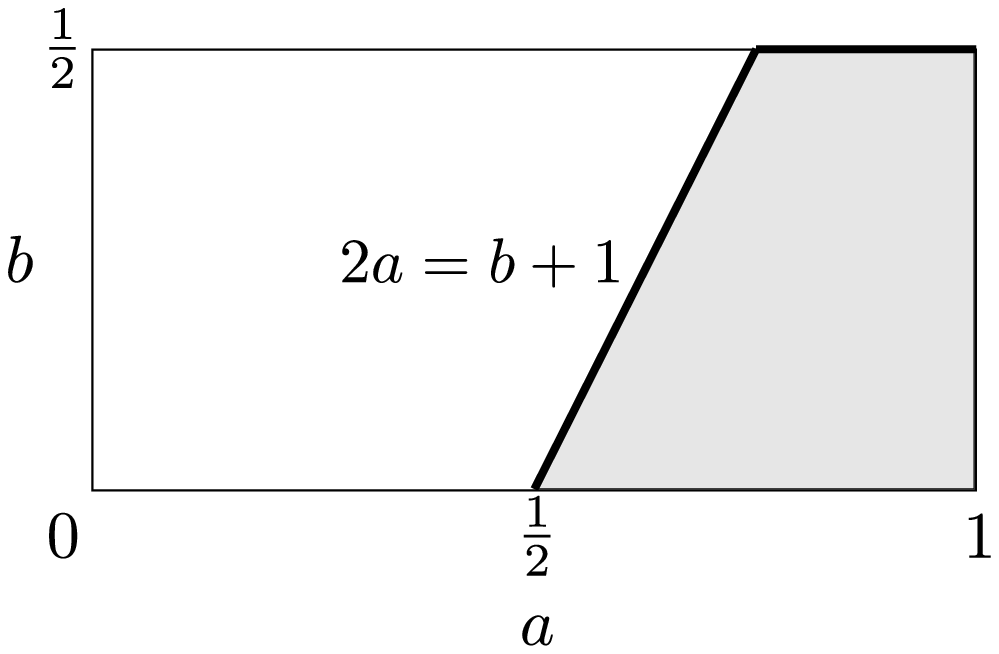}
\end{center}

\begin{definition}[Divergence Sets]\label{def:sets_Divergence}
Let $a$ and $b$ satisfy the conditions \eqref{eq:conditions_a_b}, and let $\mathcal{A}_k$, for $k\ge k_0\gg 1$, be the collection of slabs $s$ such that:
\begin{enumerate}[(i)]
\item $s$ has dimensions $cR_k^{-\frac{1}{2}}\times cR_k^{-1}\times\cdots \times cR_k^{-1}$, for $R_k = 2^k$ and $c\ll 1$.
\item $s$ has center at 
\begin{equation*}
(2p_1R_k/(qD_k^2), p_2/(D_kq),\ldots, p_n/(D_kq)),
\end{equation*}
where $(p_1/q,\ldots,p_n/q)$ is an admissible fraction (Definition~\ref{def:Admissible_fractions}) with $1\le q\le Q_k$, and $D_k$ and $Q_k$ are given by \eqref{eq:D_Q}.
\end{enumerate} 
A $(a,b)$-set of divergence $F$ is defined as
\begin{equation}\label{FSet}
F := \limsup_{k\to \infty} F_k , \qquad F_k :=  \bigcup_{s\in\mathcal{A}_k}s \, .
%\qquad 0 < c_0 \leq c \ll 1 \,.
\end{equation}
\end{definition}

For fixed $a$ and $b$, we define the initial datum
\begin{equation}\label{eq:Example}
g_{a,b} = \sum_{k\ge k_0}R_k^{-s}\frac{k}{\norm{f_{D_k}}_2}f_{D_k},
\end{equation}
where $R_k = 2^k$ and $k_0\gg 1$. Inequality \eqref{eq:thm:supremum_Fq} dictates the value of $s$, and in terms of $a$ and $b$ we have 
\begin{equation}\label{eq:s_with_a_b}
s := \frac{1}{4}+\frac{n-1}{2(n+1)}(n-(n-1)a-b).
\end{equation}
Since
\begin{equation}\label{eq:regInitdatum}
\norm{f}_{H^{s'}(\R^n)}^2 = \sum_{k\ge k_0} kR_k^{2(s'-s)}<\infty,\quad\textrm{for } s'<s,
\end{equation}
we have that $f\in H^{s'}(\R^n)$ for every $s'<s$.

We have to prove that the different terms in the sum \eqref{eq:Example} do not interfere with each other. We need the following Lemma.

\begin{lemma}\label{thm:Strong_Decay}
If the Fourier transform of $\varphi\in \mathcal{S}(\R)$ is supported in $(-1,1)$, then for every $N\ge 1$ it holds
\begin{equation}
\abs{\U \varphi(x)}\le C_N\frac{1}{\abs{x}^N},\qquad\textrm{for } \abs{x}>2t.
\end{equation}
\end{lemma}
\begin{proof}
We use the principle of non-stationary phase. We assume that $x>2t$; the other case is similar. The solution is
\begin{equation*}
\U \varphi(x) = \int \hat{\varphi}(\xi)e(-\pi t\abs{\xi}^2 + 2\pi x\xi)\,d\xi.
\end{equation*}
Since $\partial_\xi e(-\pi t\abs{\xi}^2 + 2\pi x\xi)=-2\pi i(t\xi-x)e(-\pi t\abs{\xi}^2 + 2\pi x\xi)$, then by repeated integration by parts we obtain
\begin{equation*}
\abs{\U \varphi(x)}\le C_N\frac{1}{\abs{x-t}^N},
\end{equation*}
which is the statement of the Lemma.
\end{proof}

Before proving the main result of this section, we need to make an observation on the way we define solutions. 
For $f \in H^{s}$, we define solutions for Sobolev functions, in such a way that they are well defined 
on sets with large Hausdorff dimension. Recall that $Q(N)$ is the cube of side $N$ centered at zero. We set 
\begin{equation}\label{SolSchr}
e^{it\hbar\Delta/2}f(x) = \lim_{N \to \infty} S_{N}(t)f(x) ,
\end{equation}
where
\begin{equation}\label{rt}
S_{N}(t)f(x) =
\int_{Q(N)} \hat{f}(\xi)e(-\pi t\abs{\xi}^2 + 2\pi x\cdot \xi)\,d\xi.
%\sum_{\substack{ k \in \mathbb{Z}^d \\ |k_\ell| \leq N \\ \ell=1, \ldots, d} }   \hat{f}(k)\, e^{ ik\cdot x -  i  | k |^{2} t } .
\end{equation}
The limit \eqref{SolSchr} is usually taken with respect to the $L^2$ norm, but here we take all the limits pointwise at each point $x$ where they exist. 
When $f \in L^2(\R)$, it is known that the limit exists pointwise for almost every $x \in \R$ and that it coincides with the $L^2$--limit. When $n = 1$, 
this result is due to Carleson \cite{Carleson1966},
whose proof extends to higher dimensions as proved, for instance, in~\cite{Fefferman1971}.  
Moreover, we can show that this limit exists $\gamma$-almost everywhere for every $f \in H^{s}$ with $s \in (0, n/2]$, as long 
as $\gamma>n-2s$; see the appendix of \cite{ecei2020}.  
% see Proposition \ref{Carls>0L2}. 
This can be regarded as a refinement of Carleson's result, although it does not recover it. 
%This refines Carleson's theorem, but of course it does not recover it.
%We also point out that solutions defined as in \eqref{SolSchr} coincide with the usual $L^2$--limit solutions, almost everywhere with respect to Lebesgue measure.

\begin{theorem}\label{thm:Example}
If $g_{a,b}$ is the initial datum defined in \eqref{eq:Example}, then
\begin{equation}\label{FollowsBy}
\limsup_{t \to 0^+}\abs{\U g_{a,b}(x)} = \infty
\end{equation}
for every $x\in (F\cap ( [c_0,c_1]\times [0, c_1]^{n-1})) \setminus \Omega$, where 
\begin{itemize}
\item
$c_0 := \frac{1}{10}c_1, \qquad c_1\ll c \ll 1$;
\item
$F$ is a $(a,b)$-set of divergence;
\item
$\mathcal{H}^{\gamma}(\Omega) = 0$ for $\gamma > n-2s$.
\end{itemize}
\end{theorem}

\begin{figure}[t]
\centering
\includegraphics[scale=0.7]{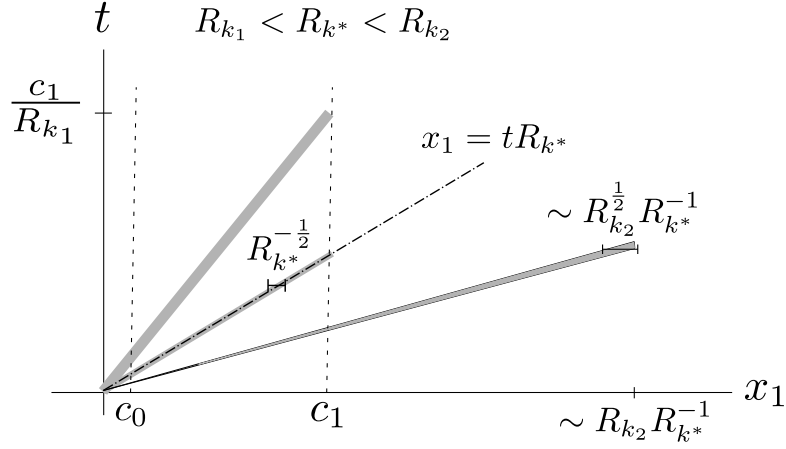}
\caption{The gray lines represent the regions where the functions $\U h_k$ concentrate.}
\end{figure}

\begin{proof}
We define $h_k := kR_k^{-s}f_{D_k}/\norm{f_{D_k}}_2$, where $R_k := 2^k$. From the proof of Theorem~\ref{thm:Supremum_Fq} we know that for 
$t = 2p_1/(D_k^2q) \ll 1/ R_k$ the value of the solution at $x\in F_k\cap [0,c_1]^n = \bigcup_{1\le q\le Q_k} E_{q,D_k}\cap [0,c_1]^n$ is
\begin{equation}\label{BiggerThan1}
\abs{\U h_k(x)} \gtrsim k.
\end{equation}
We fix $k^*\ge k_0\gg 1$ and $x\in F_{k^*}$, and we know that 
\begin{equation}\label{NonInterpherence}
x_1 = tR_{k^*} + \BigO(R_{k^*}^{-\frac{1}{2}}), \qquad c_0R_{k^*}^{-1}<t<c_1R_{k^*}^{-1}.
\end{equation} 
It suffices to prove 
\begin{equation}\label{NoInterpherence}
\abs{\U h_k(x)}\lesssim R_k^{-1}, \quad \mbox{for} \quad k\neq k^*\, ,
\end{equation}
because then for $t = 2p_1 /(D_{k^*}^2q) \ll 1/ R_{k^*}$ we would have, for all $k_1 \geq k^*$, the following (recall \eqref{rt})
\begin{equation}\label{FinalLB}
\abs{S_{2^{k_1}}(t) g_{a,b}(x)}>\abs{\U h_{k^*}(x)}-\sum_{k_0\le k\neq k^* \leq k_1}\abs{\U h_k(x)}\gtrsim k^*,
\end{equation}
as long as $k_0\gg 1$; then in order to deduce \eqref{FollowsBy} we note that for all $x \in F\cap ([c_0,c_1]\times [0, c_1]^{n-1})$ we can choose  
any $k^* \geq k_0 \gg1$, and we  
we have a lower bound as \eqref{FinalLB} and the sequence of times $t = 2p_1/(D_{k^*}^2q) \ll 1/R_{k^*}$
goes to zero as $k^* \to \infty$. More precisely, since we have 
\begin{equation}\label{FinalLB2}
e^{it\hbar\Delta/2}f(x) = \lim_{N \to \infty} S_{N}(t) g_{a,b}(x) =  \lim_{k_1 \to \infty} S_{2^{k_1}}(t) g_{a,b}(x)
\end{equation}
except possibly on sets $\Omega_{t}$, $t = 2p_1 /(D_{k^*}^2q)$ with $\mathcal{H}^{\gamma}(\Omega_t) =0$ and these sets
are countably many, then \eqref{FollowsBy} would follow by 
\eqref{FinalLB}-\eqref{FinalLB2},
taking $$\Omega := \bigcup_{t = 2p_1 /(D_{k^*}^2q)} \Omega_t.$$ 

%\begin{equation}\label{FinalLB}
%\abs{\U g_{a,b}(x)}>\abs{\U h_{k^*}(x)}-\sum_{k_0\le k\neq k^*}\abs{\U h_k(x)}\gtrsim k^*,
%\end{equation}

It remains to prove \eqref{NoInterpherence}. From \eqref{eq:f_D}, we see that we 
can bound $\U \tilde{h}_k(\tilde x)$ with the crude estimate
\begin{equation*}
\abs{\U \tilde{f}_{D_k}(\tilde x)}\lesssim \Big(\frac{R_k}{D_k}\Big)^{n-1},
\end{equation*}
so we can control each term $\U h_k(x)$ as
\begin{align*}
\abs{\U h_k(x)} &\le \abs{[\U[tR_k] \varphi(R_k^\frac{1}{2}(x_1-tR_k))]}kR_k^\frac{1}{4}\Big(\frac{R_k}{D_k}\Big)^\frac{n-1}{2}R_k^{-s} \\
&\lesssim \abs{\U[tR_k] \varphi(R_k^\frac{1}{2}(x_1-tR_k))}kQ_k^\frac{n-1}{2},
\end{align*}
and we can apply Lemma~\ref{thm:Strong_Decay} to $\varphi$.

We verify the hypotheses of Lemma~\ref{thm:Strong_Decay} when $R_k<R_{k^*}$. 
By \eqref{NonInterpherence} we get 
\begin{align*}
\frac{R^\frac{1}{2}_k(x_1-tR_k)}{tR_k} &= \frac{R_k^\frac{1}{2}(t(R_{k^*}-R_k)+\BigO(R_{k^*}^{-\frac{1}{2}}))}{tR_k} \\
&\gtrsim R_k^{-\frac{1}{2}}R_{k^*} >2, 
\\
\end{align*}
for $k,k^*\ge k_0\gg 1$; hence, $\abs{\U h_k(x)}\lesssim_N kQ_k^\frac{n-1}{2} R_k^{-\frac{N}{2}}\lesssim R_k^{-1}$, for $N\gg 1$.

We verify now the hypotheses of Lemma~\ref{thm:Strong_Decay} when $R_k>R_{k^*}$:
\begin{align*}
\frac{R^\frac{1}{2}_k(tR_k-x_1)}{tR_k} &= \frac{R_k^\frac{1}{2}(t(R_k-R_{k^*})+\BigO(R_{k^*}^{-\frac{1}{2}}))}{tR_k} \\
&\gtrsim R_k^\frac{1}{2} > 2,
\end{align*}
for $k,k^*\ge k_0\gg 1$; hence, $\abs{\U h_k(x)}\lesssim_N kQ_k^\frac{n-1}{2}(R_{k^*}R_k^{-\frac{3}{2}})^N \lesssim R_k^{-1}$, for $N\gg 1$.
\end{proof}

\section{Dimension of the Divergence Set}

In the previous section we constructed initial data parameterized by $a$ and $b$. To simplify matters, we choose those values of $a$ and $b$ for which computations are easier and exhaust all possible outcomes. Our choices are:
\begin{align}\label{Def:ab}
(\textrm{I})& &\frac{1}{2}<a\le \frac{3}{4}& &\mbox{and}& &b &= 2a - 1 \\ \label{Def:ab2}
(\textrm{II})& &\frac{3}{4}<a\le 1& &\mbox{and}& &b &= \frac{1}{2}.
\end{align}
We refer to these $(a,b)$-sets of divergence (Definition~\ref{def:sets_Divergence}) as of type I and type II. We remark that for I we have $Q = 1$, and that $a = 1$ and $b = \frac{1}{2}$ is Bourgain's example.

\begin{theorem}
Let $0 < c_0 \leq 1$. 
If $F = \limsup_{k\to\infty} F_k$ is a $(a,b)$-set of divergence (Definition~\ref{def:sets_Divergence}), then $\dim (F \cap [0,c_0]^n) \le \alpha := \frac{1}{2}+(n-1)a+b$.
\end{theorem}
\begin{proof}
Fix a scale $0<\lambda\ll 1$ and choose $k'$ such that $R_{k'}^{-1}<\lambda$. Since $F_k$ is union of $\lesssim R_k^{(n-1)a+b}$ slabs with dimensions $R_k^{-\frac{1}{2}}\times R_k^{-1}\times\cdots\times R_k^{-1}$, and each slab can be covered by $R^\frac{1}{2}_k$ balls $B_r$, for $r = R_k^{-1}$, then we can find a collection $\mathcal{B}_k$ with $|\mathcal{B}_k| = R_k^{\alpha}$ of balls with 
radius $R_k^{-1}$ covering $F_k$, so that
\begin{equation*}
\mathcal{H}^\beta_\lambda(F) := \inf\{\sum_{B_\rho\in\mathcal{B}}\rho^\beta \mid F\subset\bigcup_{B_\rho\subset\mathcal{B}}B_\rho\textrm{ and } \rho<\lambda\} \le \sum_{k\ge k'}\sum_{B_r\in\mathcal{B}_k}R_k^{-\beta},
\end{equation*}
and the last sum is smaller than $\sum_{k\ge k'}R_k^{\alpha-\beta}$, which tends to zero as $k'\to\infty$ whenever $\beta>\alpha$. 
\end{proof}

To prove the corresponding lower bound of $\dim F$, we employ the techniques in Section~4 of \cite{LucaRogers2017}. We recall a result of Falconer, which is consequence of Theorem~3.2 and Corollary~4.2 in \cite{zbMATH03979425}.

\begin{lemma}\label{thm:Falconer}
Let $0 < c \leq 1$. Suppose that there exists a constant $C >0$ such that, for all $\delta>0$ and all cubes $Q(x,\delta)\subset [0, c]^n$, we have the density condition
\begin{equation*}
\liminf_{k\to\infty}\mathcal{H}^\beta_\infty(F_k\cap Q(x,\delta))\ge C \delta^\beta,
\end{equation*}
where $\{F_k\}_{k\ge 0}$ is a sequence of open subsets of $B(0,1)$. Then, for all $\beta'<\beta$,
\begin{equation*}
\mathcal{H}^{\beta'}(\limsup_{k\to\infty}F_k)>0.
\end{equation*}
\end{lemma}

We prove now the lower bound of $\dim F$ in the easier case, in the case of sets of type I.

\begin{theorem}\label{thm:dim_TypeI}
If $F = \limsup_{k\to\infty}F_k$ is a set of type I, that is, $\frac{1}{2}<a\le\frac{3}{4}$ and $b = 2a-1$, then $\dim F \cap [0,c_0]^n \ge \alpha$ where
\begin{equation}\label{Def:alpha}
\alpha :=  \frac{1}{2}+(n-1)a+b.
\end{equation}
\end{theorem} 
\begin{proof}
From Lemma~\ref{thm:Falconer} it will be sufficient to show that
\begin{equation}\label{eq:H_Lower_bound_I}
\mathcal{H}^\beta_\infty(F_k\cap Q(x,\delta)) \geq C \delta^\beta,
\qquad \forall   Q(x,\delta) \subseteq [0,c_0]^n,
\end{equation}
holds for all $k$ sufficiently large,
where $\beta = \alpha-\varepsilon$ for $0<\varepsilon\ll 1$. The size of $k$ for which \eqref{eq:H_Lower_bound_I} holds 
will depend on $\delta$. 
To prove \eqref{eq:H_Lower_bound_I} we define an auxiliary measure which is a uniform mass measure over
$F_k\cap Q(x,\delta)$, namely 
\begin{equation*}
\mu_k(A) := \frac{\abs{A\cap F_k\cap Q(x,\delta)}}{\abs{F_k\cap Q(x,\delta)}}.
\end{equation*}
Note that $\mu_k$ depends on the set $F_k\cap Q(x,\delta)$, but we will only stress the dependence on $k$ in the notation.

Assume we have proved 
\begin{equation}\label{MainMeasComp} 
\mu_k(B_r)\le C r^\beta \delta^{-\beta}
\end{equation} 
for all sufficiently large $k$ (the size of $k$ will depend on $\delta$). 
Using \eqref{MainMeasComp} we can prove \eqref{eq:H_Lower_bound_I} easily, noting that
if $\mathcal{B}$ is a collection of balls $B_r$ that covers $F_k$, then
\begin{equation*}
1 = \mu_k(F_k\cap Q(x,\delta))\le \sum_{B_r\in\mathcal{B}}\mu_k(B_r)\le C\delta^{-\beta}\sum_{B_r\in\mathcal{B}}r^\beta.
\end{equation*}

Thus we have reduced to prove \eqref{MainMeasComp}.
To do so we have to work at several scales.
It will be useful to keep in mind that 
if $k\gg 1$ then 
\begin{equation}\label{MeasUseful} 
\abs{F_k\cap Q(x,\delta)}\simeq R_k^{\alpha - n}\delta^n
\end{equation}
and that $R_k \to \infty$ as $k \to \infty$. Many estimates below will be indeed justified taking $k$ large enough, depending on $\delta$.
\begin{enumerate}
\item Scale $r<R_k^{-1}$ 
In the worst case a ball is entirely contained in a slab from $F_k$, so
\begin{equation*}
\mu_k(B_r)\lesssim r^nR_k^{-\alpha + n}\delta^{-n} \le r^\alpha\delta^{-n} = r^\beta\delta^{-\beta} \, r^{\alpha-\beta}\delta^{\beta-n}
<r^\beta\delta^{-\beta} \, R_k^{-(\alpha-\beta)}\delta^{\beta-n};
\end{equation*}
since $\alpha-\beta>0$ and $r<R^{-1}_k$ we have 
$R_k^{-(\alpha-\beta)}\delta^{\beta-n}<1$ for $k\gg_\delta 1$ thus~\eqref{MainMeasComp} holds at this scale.
\item
Scale $R_k^{-1}<r<R_k^{-a}$.
Recall that $R^{-a}_k<R_k^{-\frac{1}{2}}$, so a ball $B_r$ cannot contain a slab. On the other hand, 
since $r<R_k^{-a}$ a ball $B_r$ intersects at most one slab, so
\begin{equation*}
\mu_k(B_r)\lesssim rR^{-(n-1)}_kR_k^{-\alpha+n}\delta^{-n} = r R_k^{-\alpha +1}\delta^{-n} = r^{\alpha}  \frac{R_k^{-(\alpha -1)}}{r^{\alpha-1}}\delta^{-n} 
  < r^\alpha\delta^{-n},
\end{equation*}
using $R_k^{-1}<r$ and $\alpha >1$. Using also $r<R^{-a}_k$ we see that 
$$\mu_k(B_r)\lesssim r^\beta R_k^{a(\beta-\alpha)}\delta^{-n}< r^\beta \delta^{-\beta}, \qquad k\gg_{\delta} 1.$$
\item Scale $R^{-a}_k<r<R^{-\frac{1}{2}}_k$. 
A ball $B_r$ intersects $\lesssim R_k^{(n-1)a}r^{n-1}$ slabs, so 
\begin{equation*}
\mu_k(B_r)\lesssim r^n R^{(n-1)a-n+1}_kR^{-\alpha+n}_k\delta^{-n}\le r^nR_k^{-b+\frac{1}{2}}\delta^{-n}.
\end{equation*}
where we used \eqref{Def:alpha}.
Since $r<R^{-\frac{1}{2}}_k$ we have that 
$$\mu_k(B_r)\lesssim r^\beta R_k^{\frac{1}{2}\beta-\frac{1}{2}n-b+\frac{1}{2}}\delta^{-n}
< r^\beta R_k^{\frac{1}{2}(\beta - \alpha)}\delta^{-n}$$ 
where we used
\begin{equation}\label{eq:proof_TypeI_A}
\alpha := (n-1)a +b +\frac{1}{2} = \frac{n-3}{2}b + \frac{n}{2}+1 + 2b -1 < n + 2b -1.
\end{equation}
Thus
$$
\mu_k(B_r)<r^\beta\delta^{-\beta}, \qquad k \gg_{\delta} 1.$$ 
 \item
Scale $R_k^{-\frac{1}{2}}<r<R_k^{-b}$.
A ball $B_r$ contains $\lesssim R_k^{(n-1)a}r^{n-1}$ slabs, so
recalling again \eqref{Def:alpha} we get
\begin{equation*}
\mu_k(B_r)\lesssim r^{n-1}R_k^{(n-1)a - n +\frac{1}{2}}R_k^{-\alpha + n}\delta^{-n} = r^{n-1}R_k^{-b}\delta^{-n}<r^{n-1+2b}\delta^{-n},
\end{equation*}
where we used $R_k^{-\frac{1}{2}}<r$. From
 $r<R_k^{-b}$ and \eqref{eq:proof_TypeI_A} we have that 
 $$
 \mu_k(B_r)\lesssim r^\beta R_k^{-b(n+2b-1-\beta)}\delta^{-n} < 
 r^\beta R_k^{-b(\alpha-\beta)}\delta^{-n}<r^\beta\delta^{-\beta}, \qquad k\gg_{\delta} 1.
 $$
\item Scale $R_k^{-b}<r<\delta$.
A ball intersects $\lesssim R_k^{(n-1)a+b}r^n$ slabs, so
\begin{equation*}
\mu_k(B_r)\lesssim r^nR_k^{(n-1)a+b-n+\frac{1}{2}}R_k^{-\alpha + n}\delta^{-n} = r^n\delta^{-n} < r^\beta\delta^{-\beta}.
\end{equation*}
\end{enumerate}
The inequality \eqref{MainMeasComp} thus holds, and so the statement of the Theorem. 
\end{proof}

The lower bound for type II sets is harder to prove, and we need a Lemma that assures us that for all $F_k$ we can find a large sub-collection of slabs uniformly distributed. Similar arguments were used in Lemma~4.3 of \cite{ecei2020} and in Sections~5.6--5.8 of \cite{Pierce2019}.

\begin{figure}
\centering
\includegraphics[scale=0.7]{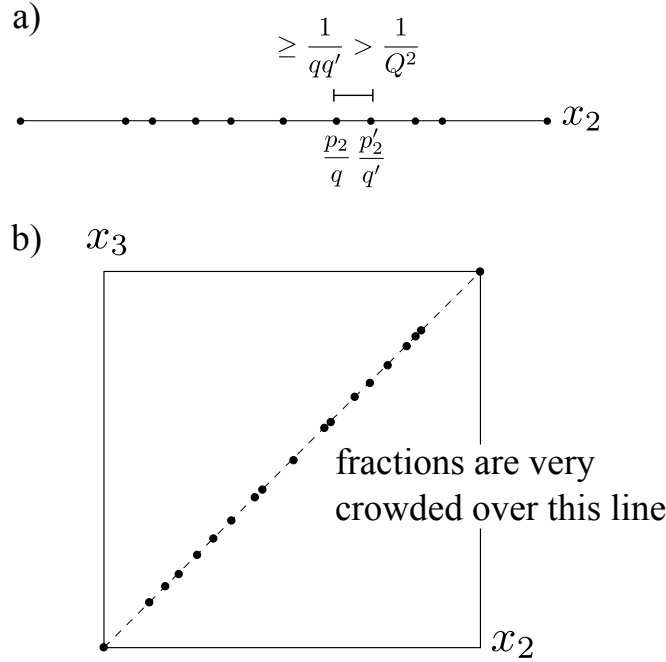}
\caption{(a) When $n = 2$ the fractions are already well separated; Lemma~\ref{thm:Selection_Slabs} is unnecessary. (b) When $n \ge 3$ the fractions might concentrate around some regions, which prohibits the Frostman measure technique we used in Lemma~\ref{thm:dim_TypeI}.}.
\end{figure} 

\begin{lemma}\label{thm:Selection_Slabs}
Let $F = \limsup_{k\to\infty}F_k$ be a set of type II, that is, $\frac{3}{4}<a\le 1$ and $b=\frac{1}{2}$. If $\mathcal{A}_k$ is the collection of slabs in $F_k\cap Q(x,\delta)$, for $\delta<1$, then, for every $\varepsilon>0$ and $k\gg_\varepsilon 1$, we can extract a sub-collection of slabs $\mathcal{A}'_k\subset\mathcal{A}_k$ such that 
\begin{enumerate}[(i)]
\item $\abs{\mathcal{A}'_k}\gtrsim R_k^{-\varepsilon}\abs{\mathcal{A}_k}$.
\item If $x = (x_1,\tilde{x})$ and $y = (y_1,\tilde{y})$ are the centers of two slabs in $\mathcal{A}'_k$ and $\tilde{x}\neq \tilde{y}$, then $\abs{\tilde{x}-\tilde{y}}\gtrsim 1/(Q^\frac{n}{n-1}_k D_k)$.
\end{enumerate} 
\end{lemma}
\begin{proof}
The sets $F_k := \bigcup_{s\in\mathcal{A}_k}s$ have a periodic structure. In fact, recall that the centers of the slabs are
\begin{equation*}
(2p_1R_k/(qD_k^2),p_2/(D_kq),\ldots, p_n/(D_kq)),
\end{equation*}
where $(p_1/q,\ldots, p_n/q)$ is an admissible fraction (Definition~\ref{def:Admissible_fractions}); hence, $F_k$ is made up of translation of the slabs in the unit cell $[0,2R_k/D_k^2]\times[0,1/D_k]^{n-1}$. We assume that $k$ is so large that the number of unit cells not entirely contained in $Q(x,\delta)$ is negligible. Therefore, the number of slabs in $Q(x,\delta)$ is $\abs{\mathcal{A}_k}\simeq D_k^{n+1}R_k^{-1}\delta^{-n}\abs{\{\textrm{slabs per unit cell}\}}$, and the Lemma reduces to extract a large number of admissible fractions in $[0,1]^{n}$ with denominator $\le Q_k$.

\begin{center}
\includegraphics[scale=0.7]{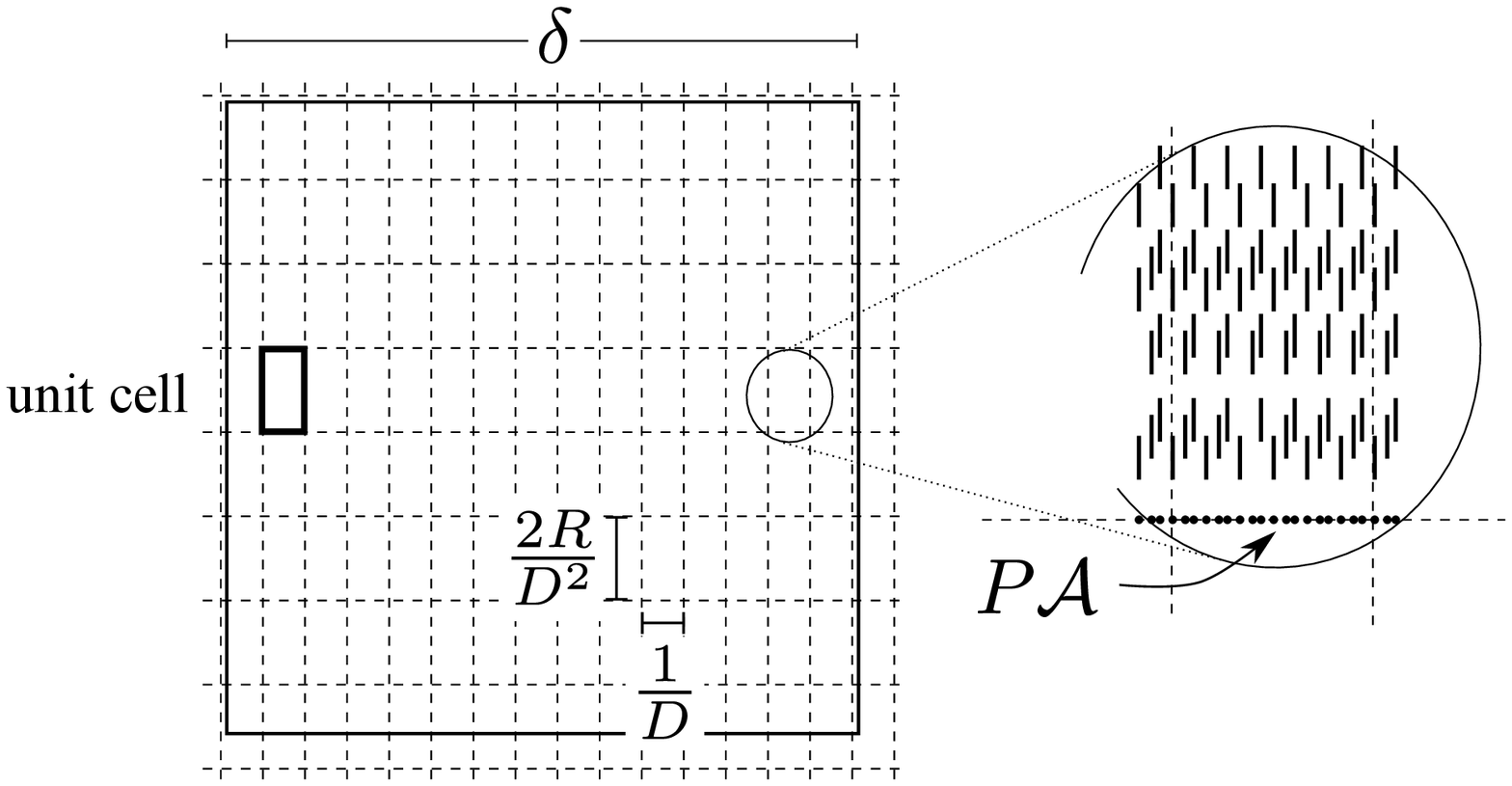}
\end{center}

We drop the subscript $k\gg 1$. Let $\mathcal{A}^0$ be the set of admissible fractions, and let $\mathcal{A}^1\subset\mathcal{A}^0$ be the collection of fractions $(p_1/q,\ldots,p_n/q)$ with $q\equiv 0\Mod 4$ and $p_j$ even for $2\le j\le n$, so that $\abs{\mathcal{A}^1}\simeq \abs{\mathcal{A}^0}$.

We denote by $P\mathcal{A}^1$ the projection of $\mathcal{A}^1$ into the plane $(x_2,\ldots,x_n)$, so $P\mathcal{A}^1$ is the set of fractions $(p_2/q,\ldots,p_n/q)$ with $q\equiv 0\Mod 4$ and even $p_j$. The Dirichlet's approximation Theorem asserts that for $2y\in\R^{n-1}$ there exists $(p_2',\ldots,p_n')\in\Int^{n-1}$ such that
\begin{equation}\label{eq:DAThm}
\abs{2y - \frac{p_j'}{q'}}\le \frac{1}{q'(Q/4)^\frac{1}{n-1}},\qquad\textrm{for some } 1\le q' \le Q/4,
\end{equation}
so if we write $q = 4q'$ and $p_j = 2p_j'$, then we can assert that for every $y\in\R^{n-1}$ there exists a fraction $(p_2/q,\ldots,p_n/q)$, for $q\equiv 0\Mod 4$ and $p_j$ even, such that
\begin{equation*}
\abs{y - \frac{p_j}{q}}\le 2^\frac{n+1}{n-1}\frac{1}{qQ^\frac{1}{n-1}},\qquad\textrm{for some } 1\le q \le Q.
\end{equation*}

In general, a point $y\in [0,1]^{n-1}$ cannot be sufficiently well approximated by fractions if it satisfies \eqref{eq:DAThm} with a fraction $(p_2'/q',\ldots,p_n'/q')$ with small $q'$, so it is convenient to ignore those points. The volume in $[0,1]^{n-1}$ occupied by those undesirable points is less than
\begin{equation}\label{eq:Undesirable_points}
\sum_{1\le q'\le Q/2^{n+2}}\Big(\frac{1}{q'(Q/4)^\frac{1}{n-1}}\Big)^{n-1}(2q')^{n-1} = \frac{1}{2}.
\end{equation}

Let $G:= \{y\in[0,1]^{n-1}\mid y\textrm{ satisfies }\eqref{eq:DAThm} \textrm{ for some } Q/2^{n+2}<q'\le Q/4\}$, then by \eqref{eq:Undesirable_points} the volume of $G$ is $>\frac{1}{2}$. Cover $G$ with cubes $Q(y,l)$, where $y\in G$ and $l := 2^{n+2+\frac{2}{n-1}}/Q^\frac{n}{n-1}$. By Vitali's covering Theorem we can find a disjoint collection of cubes $\{Q(y_j,l)\}_{1\le j\le N}$ such that
\begin{equation*}
G \subset \bigcup_{j=1}^N Q(y_j,3l);
\end{equation*}
hence, $N\ge c_nQ^n$. We pick from within each $Q(y_j,l)$ a fraction and construct so a collection of fractions $\mathcal{C}\subset P\mathcal{A}^1$; we define $\mathcal{A}^2\subset\mathcal{A}^1$ as the set of fractions such that $P\mathcal{A}^2 = \mathcal{C}$. By construction, $\abs{P\mathcal{A}^2}\gtrsim Q^n$ and any two points in $P\mathcal{A}^2$ lie at distance $\gtrsim 1/Q^\frac{n}{n-1}$; the latter, after dilation by $1/D$, implies the condition \textit{(ii)}.

The fractions in $\mathcal{A}^2$ that lie over $(p_2/q,\ldots,p_n/q)\in P\mathcal{A}^2$ is in number at least $\varphi(q)$, where $\varphi$ is the Euler's totient function. Since $\varphi(q)\ge q^{1-\varepsilon}$ for every $\varepsilon>0$ and $q\gg_\varepsilon 1$---see Theorem~327 in \cite{zbMATH05309455}---then the number of fractions in $\mathcal{A}^2$ is $\ge Q^{1-\varepsilon}\abs{P\mathcal{A}^2}\gtrsim Q^{n+1-\varepsilon} \simeq Q^{-\varepsilon}\abs{\mathcal{A}^0}$, where $\mathcal{A}^0$ is the set of admissible fractions; this concludes the verification of condition \textit{(i)}.
\end{proof}

\begin{theorem}\label{thm:dim_TypeII}
Let $0 < c_0\leq 1$. If $F = \limsup_{k\to\infty}F_k$ is a set of type II, that is, $\frac{3}{4}<a\le 1$ and $b = \frac{1}{2}$, then 
$\dim F \cap [0,c_0]^n \ge \alpha$ where
\begin{equation}\label{def:alpha2} 
\alpha := 1+(n-1)a.
\end{equation}
\end{theorem} 
\begin{proof}
We use the same method as in Theorem~\ref{thm:dim_TypeI}. For fixed $\varepsilon>0$, let $\mathcal{A}_k'$ be the collection of slabs provided 
by Lemma~\ref{thm:Selection_Slabs}, and let $F_k'$ be the corresponding set. Given $Q(x, \delta) \subseteq [0,c_0]^n$, we define again 
a measure $\mu_k$ on $F_k\cap Q(x,\delta)$ that will be useful in the proof; the measure is
\begin{equation*}
\mu_k(A) := \frac{\abs{A\cap F'_k\cap Q(x,\delta)}}{\abs{F'_k\cap Q(x,\delta)}}.
\end{equation*}
If $k\gg_\varepsilon 1$ then 
$$\abs{F_k'\cap Q(x,\delta)}\gtrsim R_k^{\alpha - n -\varepsilon}\delta^n.$$
 We take $\beta := \alpha-2n\varepsilon <\alpha - \varepsilon$.
The goal is again to prove \eqref{MainMeasComp}, from 
which we deduce Theorem \ref{thm:dim_TypeII} proceeding as we did in the proof of Theorem \ref{thm:dim_TypeI}.

Since $b=\frac{1}{2}$, we can think of the slabs over $(p_2/(qD_k),\ldots,p_n/(qD_k))$ as a single tube of length 1.
\begin{enumerate}
\item Scale $r<R^{-1}$. 
In the worst case a ball is entirely contained in a slab from $F_k$, so
\begin{equation*}
\mu_k(B_r)\lesssim r^nR_k^{-\alpha + n +\varepsilon}\delta^{-n} \le r^\alpha\delta^{-n} = r^\beta\delta^{-\beta}(r^{\alpha-\beta-\varepsilon}\delta^{\beta-n});
\end{equation*}
since $\alpha-\beta>\varepsilon$ and $r<R^{-1}_k$, then $\mu_k(B_r)<r^\beta\delta^{-\beta}$ whenever $k\gg_{\delta} 1$.
\item
Scale $R_k^{-1}<r<R_k^{-a}$.
By the properties of separation of the slabs in $\mathcal{A}_k'$, a ball $B_r$ intersects at most one slab---recall Lemma~\ref{thm:Selection_Slabs}(ii) and \eqref{eq:Condition_a}---so
\begin{equation*}
\mu_k(B_r)\lesssim rR_k^{-(n-1)-\alpha + n +\varepsilon}\delta^{-n} = rR_k^{-\alpha + 1 + \varepsilon}\delta^{-n}< r^{\alpha-\varepsilon}\delta^{-n},
\end{equation*}
where we used $\alpha >1$.
Since $r<R^{-a}_k$ we see that
 $$\mu_k(B_r)\lesssim r^\beta R_k^{a(\beta-\alpha+\varepsilon)}\delta^{-n}, \qquad  k \gg_{\delta} 1.$$
\item
Scale $R_k^{-a}<r<R_k/D_k^2 = R_k^{\frac{n-1}{n+1}(2a-\frac{3}{2})-\frac{1}{2}}$. 
A ball intersects $\lesssim R_k^{(n-1)a}r^{n-1}$ ``tubes'' of length 1 and radius $R_k^{-1}$, so (recall \eqref{def:alpha2})
\begin{equation*}
\mu_k(B_r)\lesssim r^nR_k^{(n-1)a-(n-1)-\alpha+n+\varepsilon}\delta^{-n} = r^nR_k^\varepsilon\delta^{-n} = r^\beta\delta^{-\beta}(r^{n-\beta}R_k^\varepsilon\delta^{\beta-n});
\end{equation*}
since $r<R_k/D_k^2\le R_k^{-\frac{1}{n+1}}$, we see that 
$$\mu_k(B_r)\lesssim r^\beta\delta^{-\beta}, \qquad k \gg_{\delta} 1.
$$
\item
Scale $R_k/D_k^2<r<\delta$. 
A ball $B_r$ contains $\simeq D_k^{n+1}R^{-1}_kr^n$ translations of the unit cell $[0,2R_k/D_k^2]\times [0,1/D_k]^{n-1}$. If $V$ is the volume of $F_k'$ per unit cell, then $\abs{B_r\cap F_k'}\simeq VD_k^{n+1}R^{-1}_kr^n$ and 
$$\abs{Q(x,\delta)\cap F_k'}\simeq VD_k^{n+1}R^{-1}_k\delta^n;$$
hence
\begin{equation*}
\mu_k(B_r)\lesssim r^n\delta^{-n}<r^\beta\delta^{-\beta}.
\end{equation*}
\end{enumerate}
The inequality $\mu_k(B_r)\le Cr^\beta\delta^{-\beta}$ holds for $k$ sufficiently large (depending on~$\delta$), so the proof is complete.
\end{proof}

\section{Conclusion of the proof}

We are now ready to prove our statement combining the results from the previous section. 
First we take $a, b$ as in \eqref{Def:ab}-\eqref{Def:ab2} and recall that we have defined
\begin{equation}\label{FinalChoiceOfAlpha}
\alpha := \frac{1}{2}+(n-1)a+b.
\end{equation}
Note that we have a bijection between $a \in (1/2, 1]$ (which predicts also the value of $b$ by \eqref{Def:ab}-\eqref{Def:ab2}) and 
 $\alpha \in (n/2, n]$, which is the range we are interested in (the case $\alpha = n$ was handled in \cite{Bourgain2016}). 
%Moreover, \eqref{FinalChoiceOfAlpha} reduces to
% \begin{equation}\label{RealFinalChoiceOfAlpha}
% \alpha := \left\{
% \begin{array}{lll}
% (n+1) a - \frac12 & \mbox{for} & a \in \left(\frac12, \frac34 \right]
% \\
% &&
% \\
%  (n-1) a +1 & \mbox{for} & a \in \left(\frac34, 1 \right] \,.
%    \end{array}\right.
% \end{equation} 

First we claim that given any
 \begin{equation}\label{UDC1}
  s' < s := \frac{n}{2(n+1)}+\frac{n-1}{2(n+1)}(n-\alpha)
 \end{equation}
we can find a solution $u(x,t)$ with initial datum $u_0 \in H^{s'}(\R^n)$ such that 
$$\limsup_{t \to 0^+} |u(x,t)| = \infty$$
for $ x\in (F\cap ( [c_0,c_1]\times [0, c_1]^{n-1})) \setminus \Omega$, where $F$ is an $(a,b)$-set of divergence,
$0 < c_0 := \frac{1}{10} c_1 \ll 1$ and $\Omega$ has dimension $\leq n-2s$. 
Indeed, it suffices to choose $u_0 := g_{a,b}$ defined in 
 \eqref{eq:Example} so that $u_0\in H^{s'}(\R^n)$ for 
 \begin{equation}\label{UDC2}
s' < s := \frac{1}{4}+\frac{n-1}{2(n+1)}(n-(n-1)a-b) ;
\end{equation}
see \eqref{eq:regInitdatum}-\eqref{eq:s_with_a_b}. 
Since  under \eqref{FinalChoiceOfAlpha} 
the inequality 
\eqref{UDC2} becomes \eqref{UDC1}, then the claim follows invoking 
Theorem \ref{thm:Example}. 

 Thus, to conclude the proof, we need to show that
\begin{equation}\label{VeryFinal}
\dim \left(  (F\cap ( [c_0,c_1]\times [0, c_1]^{n-1})) \setminus \Omega \right) \geq \alpha.
\end{equation}
First, covering $(F\cap ( [c_0,c_1]\times [0, c_1]^{n-1}))$ with $\simeq (c_1/c_0)^{n-1}$ cubes of side $c_0$, we see as consequence 
of Theorems \ref{thm:dim_TypeI} and \ref{thm:dim_TypeII} that
$$
\dim  (F\cap ( [c_0,c_1]\times [0, c_1]^{n-1}))  \geq \alpha.
$$
On the other hand, we know that $\dim \Omega \leq n-2s$ (see Theorem \ref{thm:Example}). Thus, since for our choice \eqref{UDC1} of $s$ we have
$\alpha > n-2s$ when $\alpha > n/2$, then \eqref{VeryFinal} follows and the proof is concluded.   

%
%If we replace $\alpha := \dim F = (n-1)a+b+\frac{1}{2}$ in \eqref{eq:s_with_a_b} for the range of type I sets and use Theorem~\ref{thm:Example}, then we get Theorem~\ref{thm:Main} for $\frac{n+1}{8}\le s<\frac{n}{4}$.

%If we replace $\alpha := \dim F = (n-1)a+b+\frac{1}{2}$ in \eqref{eq:s_with_a_b} for the range of type II sets and use Theorem~\ref{thm:Example}, then we get Theorem~\ref{thm:Main} for $\frac{n}{2(n+1)}\le s<\frac{n+1}{8}$, and we recover so the results in \cite{zbMATH07036806}.

\bibliographystyle{plain}
\bibliography{Conv_references}

\end{document}